\documentclass{article}
\usepackage{amssymb, amsmath, amsfonts}

\begin{document}

\begin{center}
{\Large Some remarks regarding special elements in algebras obtained by the
Cayley-Dickson process over Z}$_{p}$%
\begin{equation*}
\end{equation*}

Cristina Flaut and Andreea Baias%
\begin{equation*}
\end{equation*}
\end{center}

\textbf{Abstract.} {\small In this paper we provide some properties of }$k$%
{\small -potent elements in algebras obtained by the Cayley-Dickson process
over }$\mathbb{Z}_{p}${\small . Moreover, we find a structure of nonunitary
ring over Fibonacci quaternions over \thinspace }$\mathbb{Z}_{3}${\small \
and we present a method to encrypt plain texts, by using invertible elements
in such algebras.}%
\begin{equation*}
\end{equation*}

\textbf{1. Preliminaries}

\begin{equation*}
\end{equation*}%
\qquad \qquad\ \ \ 

In [MS; 11], the authors provided some properties regarding quaternions over
the field $\mathbb{Z}_{p}$. Since quaternions are special cases of algebras
obtained by the Cayley-Dickson process, in this paper we extend the study of 
$k$-potent elements over quaternions to an arbitrary algebra obtained by the
Cayley-Dickson process. These algebras, in general, are poor in properties:
are not commutative, starting with dimension $4$ (the quaternions), are not
associative, strating with dimension $8$ (the octonions) and lost
alternativity, starting with dimension $16$ (the sedionions).The good news
is that all algebras obtained by the Cayley-Dickson process are
power-associative and this is the property which will be used when we study
the $k$-potent elements in these algebras. The paper is organized as
follows: in Introduction, we present basic properties of algebras obtained
by the Cayley-Dickson process, in section 3, we characterize the $k$-potent
elements in these algebras, in section 4, we give a structure of non-unitary
and noncommutative ring over the Fibonacci quaternions over $\mathbb{Z}_{3}$
and in the last section, we provide an encryption method by using invertible
elements from these algebras.

\bigskip

\begin{equation*}
\end{equation*}

\textbf{2. Introduction}%
\begin{equation*}
\end{equation*}

\bigskip

In the following, we consider $A,$ a finite dimensional unitary algebra over
a field $\ K$ with $charK\neq 2$.

An algebra $A$ is called \textit{alternative} if $x^{2}y=x\left( xy\right) $
and $xy^{2}=\left( xy\right) y,$ for all $x,y\in A,$ \textit{\ flexible} if $%
x\left( yx\right) =\left( xy\right) x=xyx,$ for all $x,y\in A$ and \textit{%
power associative} if the subalgebra $<x>$ of $A$ generated by any element $%
x\in A$ is associative. $\ $Each alternative algebra is$\ $a\ flexible
algebra and a power associative algebra.

We consider the algebra $A\neq K$ such that for each element $x\in A$, the
following relation is true 
\begin{equation*}
x^{2}+t_{x}x+n_{x}=0,
\end{equation*}%
for all $x\in A$ and $t_{x},n_{x}\in K$. This algebra is called a \textit{%
quadratic algebra}.

It is well known that a finite-dimensional algebra $A$ is \textit{a division}
algebra if and only if $A$ does not contain zero divisors (See \textbf{[}%
Sc;66\textbf{]}).

A \textit{composition algebra} $A$ over the field $K$ is an algebra, not
necessarily associative, with a nondegenerate quadratic form $n$ which
satisfies the relation

\begin{equation*}
n(xy)=n(x)n(y),\forall x,y\in A.
\end{equation*}

A unital composition algebras are called \textit{Hurwitz algebras}.\medskip

\textbf{Hurwitz's Theorem.[}Ba; 01\textbf{]} $\mathbb{R}$, $\mathbb{C}$, $%
\mathbb{H}$ \textit{and} $\mathbb{O}$ \textit{are the only real alternative
division algebras}.\medskip

\textbf{Theorem 1.} (Theorem 2.14, \textbf{[}McC,80\textbf{]}) $A$ \textit{%
is a composition algebra if and only if} $A$ \textit{is an alternative
quadratic algebra.\smallskip }

An element $x$ in a ring $R$ is called \textit{nilpotent }if we can find a
positive integer $n$ such that $x^{n}=0$. The number $n$, the smallest with
this property, is called the \textit{nilpotency index}. A power-associative
algebra $A$ is called a \textit{nil algebra} if and only if each element of
this algebra is nilpotent.\ An element $x$ in a ring $R$ is called $k$%
\textit{-potent}, for $k>1$, a positive integer, if $k$ is the smallest
number such that $x^{k}=x$. The number $k$ is called the $k$\textit{-potency
index}. For $k=2,$ we have idempotent elements, for $k=3$, we have tripotent
elements, etc.

Let $A$ be an algebra over the field $K$ and a \textit{scalar} \textit{%
involution} $\,$over $A,$%
\begin{equation*}
\,\,\,\overline{\phantom{x}}:A\rightarrow A,a\rightarrow \overline{a},
\end{equation*}%
$\,\,$ that means a linear map with the following properties$\,\,\,$%
\begin{equation*}
\overline{ab}=\overline{b}\overline{a},\,\overline{\overline{a}}=a,
\end{equation*}%
$\,\,$and 
\begin{equation*}
a+\overline{a},a\overline{a}\in K\cdot 1,\ \text{for all }a,b\in A.\text{ }
\end{equation*}%
For the element $a\in A$, the element $\,\overline{a}$ is called the \textit{%
conjugate} of the element $a$. The linear form$\,\,$%
\begin{equation*}
\,\,\mathbf{t}:A\rightarrow K\,,\,\,\mathbf{t}\left( a\right) =a+\overline{a}
\end{equation*}%
and the quadratic form 
\begin{equation*}
\mathbf{n}:A\rightarrow K,\,\,\mathbf{n}\left( a\right) =a\overline{a}\ 
\end{equation*}%
are called the \textit{trace} and the \textit{norm \ }of \ the element $a,$
respectively$.$ From here, it results that an algebra $A$ with a scalar
involution is a quadratic algebra. $\,$Indeed, if in the relation $\mathbf{n}%
\left( a\right) =a\overline{a}$, we replace $\overline{a}=\,\mathbf{t}\left(
a\right) -a$, we obtain 
\begin{equation}
a^{2}-\,\mathbf{t}\left( a\right) a+\mathbf{n}\left( a\right) =0.  \tag{1.}
\end{equation}

Let$\,\,\,\delta \in K$ \thinspace be\ a fixed non-zero element. We define
the following algebra multiplication on the vector space $A\oplus A$ 
\begin{equation}
\left( a_{1},a_{2}\right) \left( b_{1},b_{2}\right) =\left(
a_{1}b_{1}+\delta \overline{b}_{2}a_{2},a_{2}\overline{b_{1}}%
+b_{2}a_{1}\right) .  \tag{2.}
\end{equation}%
\newline
The obtained algebra structure over $A\oplus A,$ denoted by $\left( A,\delta
\right) ,$is called the \textit{algebra obtained from }$A$\textit{\ by the
Cayley-Dickson process.} $\,$We have that $\dim \left( A,\delta \right)
=2\dim A$.

Let $x\in \left( A,\delta \right) $, $x=\left( a_{1},a_{2}\right) $. The map 
\begin{equation*}
\,\,\,\overline{\phantom{x}}:\left( A,\delta \right) \rightarrow \left(
A,\delta \right) \,,\,\,x\rightarrow \bar{x}\,=\left( \overline{a}%
_{1},-a_{2}\right) ,
\end{equation*}%
\newline
is a scalar involution of the algebra $\left( A,\delta \right) $, extending
the involution $\overline{\phantom{x}}\,\,\,$of the algebra $A.$ We consider
the maps 
\begin{equation*}
\,\mathbf{t}\left( x\right) =\mathbf{t}(a_{1})
\end{equation*}%
and$\,\,\,$ 
\begin{equation*}
\mathbf{n}\left( x\right) =\mathbf{n}\left( a_{1}\right) -\delta \mathbf{n}%
(a_{2})
\end{equation*}%
called $\,\,$the \textit{trace} and the \textit{norm} of the element $x\in $ 
$\left( A,\delta \right) ,$ respectively.\thinspace $\,$

\thinspace If we consider $A=K$ \thinspace and we apply this process $t$
times, $t\geq 1,\,\,$we obtain an algebra over $K,\,\,$%
\begin{equation}
A_{t}=\left( \frac{\delta _{1},...,\delta _{t}}{K}\right) .  \tag{3. }
\end{equation}

Using induction in this algebra, the set $\{1,f_{1},...,f_{n-1}\},n=2^{t},$
generates a basis with the properties:%
\begin{equation}
f_{i}^{2}=\delta _{i}1,\,\,_{i}\in K,\delta _{i}\neq 0,\,\,i=1,...,t 
\tag{4.}
\end{equation}%
and \ 
\begin{equation}
f_{i}f_{j}=-f_{j}f_{i}=\alpha _{ij}f_{k},\,\,\alpha _{ij}\in K,\,\,\alpha
_{ij}\neq 0,i\neq j,i,j=\,\,1,...n-1,  \tag{5.}
\end{equation}%
$\ \ \alpha _{ij}$ and $f_{k}$ being uniquely determined by $f_{i}$ and $%
f_{j}.$

From \textbf{[}Sc; 54\textbf{]}, Lemma 4, it results that in any algebra $%
A_{t}$ with the basis \newline
$\{1,f_{1},...,f_{n-1}\}$ satisfying relations $\left( 4\right) $ and $%
\left( 5\right) ,$ we have:

\begin{equation}
f_{i}\left( f_{i}x\right) =\delta _{i}x=(xf_{i})f_{i},  \tag{6.}
\end{equation}%
for all $i\in \{1,2,...,n-1\}$ and for \ every $x\in A.$

The field $K$ is the center of the algebra $A_{t},$for $t\geq 2.$(See [Sc;
54]). Algebras $A_{t}$ of dimension $2^{t}\ $obtained by the Cayley-Dickson
process, described above, are flexible and\textit{\ }power associative for
all $t\geq 1$ and, in general, are not division algebras for all $t\geq 1$.\ 

For $t=2,$ we obtain the generalized quaternion algebras over the field $K$.
If we take $K=\mathbb{R}$ and $\delta _{1}=\delta _{2}=-1$, we obtain the
real quaternion algebra over $\mathbb{R}.$ This algebra is an associative
and a noncommutative algebra and will be denoted with $\mathbb{H}$.

Let $\mathbb{H}$ be the real quaternion algebra with basis $\{1,i,j,k\}$,
where 
\begin{equation}
i^{2}=j^{2}=k^{2}=-1,ij=-ji,ik=-ki,jk=-kj.  \tag{7.}
\end{equation}%
Therefore, each element from $\mathbb{H}$ has the following form$~\ $%
\begin{equation*}
q=a+bi+cj+dk,a,b,c,d\in \mathbb{R}.
\end{equation*}

We remark that $\mathbb{H}$ is a vector space of dimension $4$ over $\mathbb{%
R}$ with the addition and scalar multiplication. Moreover, $\mathbb{H}$ has
a ring structure with multiplication given by $\left( 7\right) $ and the
usual distributivity law.

If we consider $K$ a finite field with $charK\neq 2$, due to the
Wedderburn's Theorem, a quaternion algebra over $K$ is allways a non
division algebra or a split algebra.

\begin{equation*}
\end{equation*}

\textbf{3. Characterization of }$k$\textbf{-potent elements in algebras
obtained by the Cayley-Dickson process} 
\begin{equation*}
\end{equation*}

In the paper \textbf{[}Mo; 15\textbf{]}, the author gave several
characterizations of $k$-potent elements in associative rings from an
algebraic point of view. In \textbf{[}RPC; 22\textbf{]}, the authors
presented some properties of $\left( m,k\right) $-type elements over the
ring of integers modulo $n$ and in \textbf{[}Wu; 10\textbf{]}, the author
emphasize the applications of $k$-potent matrices to digital image
encryption.

In the following, we will study the properties of $k$-potent elements in a
special case of nonassociative structures, that means we characterize the $k$%
-potent elements in algebras obtained by the Cayley-Dickson process over the
field of integers modulo $p$, $p$ a prime number greater than $2$, $K=%
\mathbb{Z}_{p}$.\smallskip

\textbf{Remark 2.} Since algebras obtained by the Cayley-Dickson process are
power associative, we can define the power of an element. In this paper, we
consider $A_{t}$ such an algebra, given by the relation $\left( 3\right) $,
with $\delta _{i}=-1$, for all $i$, $i\in \{1,...,t\}$. We consider $x\in
A_{t},$ a $k$-potent element, that means $k$ is the smallest positive
integer with this property. Since $A_{t}$ is a quadratic algebra, from
relation $\left( 1\right) $, we have that $x^{2}-\,\mathbf{t}\left( x\right)
x+\mathbf{n}\left( x\right) \,=0$, with $\,\mathbf{t}\left( x\right) \in K$
the trace and $\mathbf{n}\left( x\right) \in K$ the norm of the element $x$.
To make calculations easier, we will denote $\mathbf{t}\left( x\right)
=t_{x} $ and $\mathbf{n}\left( x\right) =n_{x}$.\smallskip

\textbf{Remark 3.} In general, algebras obtained by the Cayley-Dickson
process are not composition algebras, but the following relation 
\begin{equation*}
\mathbf{n}\left( x^{m}\right) =\left( \mathbf{n}\left( x\right) \right) ^{m}
\end{equation*}%
is true, for \thinspace $m$ a positive integer. Indeed, we have $\mathbf{n}%
\left( x^{m}\right) =x^{m}\overline{x^{m}}$ and $\left( \mathbf{n}\left(
x\right) \right) ^{m}=(x\overline{x})^{m}=x\overline{x}\cdot ...\cdot x%
\overline{x}$, $m$-times with $\overline{x}=t_{x}-x,t_{x}\in K$. Since $x$
and$~\overline{x}$ are in the algebra generated by $x$, they associate and
comute, due to the power associativity property. If $x\in A_{t}$ is an
invertible element, that means $n_{x}\neq 0$, then the same remark is also
true for $x^{-1}=\frac{\overline{x}}{n_{x}}$, the inverse of the element $x$%
. The element $x^{-1}$ is in the algebra generated by $x$, therefore
associate and comute with $x$.

ii) We know that $x^{2}-\,t_{x}x+n_{x}\,=0$. If $x\in A_{t}$ is a nonzero $k$%
\textit{-}potent element, then, from the above, we have $n_{x}=0$ or $%
n_{x}\neq 0$ and $n_{x}^{k-1}=1$.

iii) Let $x\in A_{t}$ be a nonzero $k$\textit{-}potent element such that $%
n_{x}\neq 0$. Then, the element $x$ is an invertible element in $A_{t}$ such
that\textbf{\ }$x^{k-1}=1$. \textbf{\ }Indeed, if $x^{k}=x$, multiplying
with \textbf{\smallskip }$x^{-1}$ we have$~$ $x^{k-1}=1$.

iv) For a nilpotent element $x\in A_{t}$ there is a positive integer $k\geq
2 $ such that $x^{k}=0,k$ the smallest with this property. From here, we
have that $n_{x}=0$, therefore $x^{2}=t_{x}x$. It results that $%
x^{k}=t_{x}x^{k-1} $, then $t_{x}x^{k-1}=0$ with $x^{k-1}\neq 0$. We get
that $t_{x}=0$ and \thinspace $x^{2}=0$. Therefore, we can say that in an
algebra $A_{t}$, if exist, we have only nilpotent elements of index
two.\smallskip

In the following, we will characterize the $k$-potent elements in the case
when $n_{x}=0$.\smallskip

\textbf{\ Proposition 4. }\textit{The element} $x\in A_{t}$, $x\neq 0$, $\ $%
\textit{with} $n_{x}=0$ \textit{and} $t_{x}\neq 0$ \textit{is a} $k$\textit{%
-potent element in }$A_{t}$ \textit{if and only if} $t_{x}$ \textit{is a} $k$%
\textit{-potent element in} $\mathbb{Z}_{p}^{\ast },2\leq k\leq p$ ($t_{x}$ 
\textit{has} $k-1$ \textit{as multiplicative order in} $\mathbb{Z}_{p}^{\ast
}$).\smallskip

\textbf{Proof. }We must prove that if $k$ is the smallest positive integer
such that $x^{k}=x$, then \thinspace $t_{x}^{k}=t_{x}$, therefore \thinspace 
$t_{x}^{k-1}=1$, with $k$ the smallest positive integer with this property.

We have $%
x^{k}=x^{k-2}x^{2}=x^{k-2}t_{x}x=t_{x}x^{k-1}=t_{x}x^{k-3}x^{2}=t_{x}^{2}x^{k-2}=...=t_{x}^{k-1}x 
$. If $t_{x}^{k-1}=1$, we have $x^{k}=x$ and if $x^{k}=x$, we have $%
x=t_{x}^{k-1}x$, therefore $t_{x}^{k-1}=1$.

Now, we must prove that $k\leq p$. We know that in $\mathbb{Z}_{p}$ the
multiplicative order of a nonzero element is a divisor of $p-1$. If the
order is $p-1$, the element is called a primitive element. If $t_{x}\neq 0$
in $\mathbb{Z}_{p}$ and $t_{x}^{k-1}=1$, it results that $(k-1)\mid (p-1)$,
then $k-1\leq p-1$ and $k\leq p$.\smallskip

\textbf{Remark 5.} For\textit{\ }elements $x$ with $n_{x}=0$ and $t_{x}\neq
0 $, from the above theorem, we remark that in an algebra $A_{t}\,$\ over $%
\mathbb{Z}_{p}$\ we have $k\leq p$, where $k$\textit{\ }is the potency
index. That means the $k$-potency index in these conditions does not exceed
the prime number $p$. Since $a^{p-1}\equiv 1$ \textit{mod} $p$, for all
nonzero $a\in \mathbb{Z}_{p}$, allways it results that $x^{p}=x$. It is not
necessary for $p$ to be the smallest with this property.\smallskip\ 

\textbf{Example 6.} If we take $p=5$ and we have $x\in A_{t}$ such that $%
x^{38}=x$, since we known that $x^{5}=x$, we obtain $x^{38}=x^{35}x^{3}=%
\left( x^{5}\right) ^{7}x^{3}=x^{7}x^{3}=x^{10}=x^{5}x^{5}=x^{2}$.
Therefore, $x^{2}=x$ and the $k$-potency index is $2$.\smallskip\ 

In the following, we will characterize the $k$-potent elements when $%
n_{x}\neq 0$ and $n_{x}^{k-1}=1$. We suppose that $k\geq 3$. Indeed, if $k=2$%
, we have $x^{2}=x$, then $x=1$.\smallskip

The following result it is well known from literature. We reproduce here the
proof\textbf{.\smallskip }

\textbf{Proposition 7.} \textit{Each element of a finite field} $K$\textit{\
can be expressed as a sum of two squares from }$K$.\smallskip

\textbf{Proof.} If $charK=2$, we have that the map $f:K\rightarrow K,f\left(
x\right) =x^{2}$ is an injective map, therefore is bijective and each
element from $K$ is a square. Indeed, if $f\left( x\right) =f\left( y\right) 
$, we have that $x^{2}=y^{2}$ and $x=y$ or $x=-y=y,$since $-1=1$ in $charK=2$%
.

Assuming that $charK=p\neq 2$. We suppose that $K$ has $q=p^{n}$ elements,
then $K^{\ast }$ has $q-1$ elements. Since $\left( K^{\ast },\cdot \right) $
is a cyclic group with $q-1$ elements, \thinspace $K^{\ast
}=\{1,v,v^{2},...,v^{q-2}\}$, half of them, namely the even powers are
squares. The zero element is also a square, then we have $\frac{q-1}{2}+1=%
\frac{q+1}{2}$ square elements from \thinspace $K$ which are squares. We
known that from a finite group $\left( G,\ast \right) $ if $S$ and $T$ are
two subsents of $G$ and $\left\vert S\right\vert +\left\vert T\right\vert
>\left\vert G\right\vert $, we have that each $x\in G$ can be expresses as $%
x=s\ast t,$ $s\in S,t\in T$. For $g\in G$, we consider the set $%
gS^{-1}=\{g\ast s^{-1},$ $s\in S\}$ wich has the same cardinal as the set $T$%
. Since $\left\vert S\right\vert +\left\vert T\right\vert >\left\vert
G\right\vert $, it results that $\left\vert T\right\vert +\left\vert
gS^{-1}\right\vert >\left\vert G\right\vert $, therefore $T\cap gS^{-1}\neq
\emptyset $. Then, there are the elements $s\in S$ and $t\in T$ such that $%
t=g\ast s^{-1}$ and $g=s\ast t$. Now, if we consider $S$ and $T$ two sets
equal with the multiplicative. In the group $\left( K,+\right) $, we have
that $\left\vert S\right\vert +\left\vert T\right\vert =q+1>\left\vert
K\right\vert $, therefore each $x\in K$ can be writen as $x=s^{2}+t^{2}$,
with $s\in S,t\in T$.\smallskip

\textbf{Remark 8.} i) We can find an element $w\in A_{t}$, different from
elements of the base, such that $w^{2}=-1$. Indeed, such an element has $%
n_{w}=1$ and $t_{x}=0$. With the above notations and from the above
proposition, since $1=a^{2}+b^{2}$, we can take $w_{ij}=af_{i}+bf_{j}$, $%
a,b\in \mathbb{Z}_{p}$ and $f_{i},f_{j}$ elements from the basis in $A_{t}$,
given by $\left( 4\right) $. Therefore, $w_{ij}^{2}=-1$.

ii) The group $\left( \mathbb{Z}_{p}^{\ast },\cdot \right) $ is cyclic and
has $p-1$ elements. Elements of order $p-1$ are primitive elements. The rest
of the elements have orders divisors of $p-1$.\smallskip

Now, we consider the equation in $A_{t}$ 
\begin{equation}
x^{n}=1,n\text{ a positive integer.}  \tag{8.}
\end{equation}%
\qquad \qquad

In the following, we will find some conditions such that this equation has
solutions different from $1$.

\textbf{Remark 9. } i) With the above notations, we consider $w\in A_{t}$ a
nilpotent element (it has the norm and the trace zero). Therefore, the
element $z=1+w$ has the property that $z^{n}=1+nw$, therefore if $n=pr,r$ a
positive integer, the equation $\left( 8\right) $ has solutions of the form $%
z=1+w$, for all nilpotent elements $w\in A_{t}$. It is clear that $z$ has
the norm equal with $1$ and $z^{p}=1$, therefore $z^{p+1}=z$, is a $p$%
-potent element.

ii) If we consider $\eta \in \mathbb{Z}_{p}^{\ast }$ with the multiplicative
order $\theta $ and $z=\eta +w,\,w$ nilpotent,$\ $we have that $\left( \eta
+w\right) ^{p}=\eta ^{p}+pw=\eta ^{p}\,\ $and $\left( \eta +w\right)
^{p\theta }=1$. Therefore, if $n=pr,r$ a positive integer, the equation $%
\left( 8\right) $ has solutions of the form $z=1+w$, for all nilpotent
elements $w\in A_{t}$. If $r$ is a multiplicative order of an element from
\thinspace $\mathbb{Z}_{p}^{\ast }$ and $n=pr,r$ a positive integer, then
the equation $\left( 8\right) $ has solutions of the form $z=\eta +w$, for
all $\eta \in A_{t}$, $\eta $ of order $r,w\,\ $a nilpotent element in $%
A_{t} $.

iii)With the above notations, we consider the element $w\in A_{t}$ sucht
that $w^{2}=-1$ and $z=1+w$. We have that $z^{2}=\left( 1+w\right)
^{2}=2w,z^{3}=\left( 1+w\right) ^{3}=2w-2$ and $z^{4}=\left( z^{2}\right)
^{2}=-4$ \textit{modulo} $p$. Let $\eta =-4\in \mathbb{Z}_{p}^{\ast }$ with
the multiplicative order $\theta $, $\theta $ is allways \qquad \qquad
\qquad an even number. We have that $z^{4\theta }=1$.

iv) Let \textit{\ \thinspace }$z=a+w\in A_{t}$, where\textit{\ }$a\in 
\mathbb{Z}_{p}$\textit{\ }and\textit{\ }$w\in A_{t}$\textit{, }with $t_{w}=0$
and $n_{w}\neq 0$. We have that $w^{2}=\alpha $\textit{\ }$\in \mathbb{Z}%
_{p}^{2}$, therefore, $z^{r}=C_{r}+D_{r}w$. If $z^{s}=1$, then there is a
positive integer $m\leq s$ such that $C_{m}=0$ or $D_{m}=0$. Indeed, if $m=s$%
, we have $D_{s}=0$ and $C_{s}=1$.\smallskip

\textbf{Proposition 10.} \textit{By using the above notations, we consider
the element \thinspace }$z=a+w$\textit{, where }$a\in \mathbb{Z}_{p}$\textit{%
\ and }$w\in A_{t}$\textit{\ with the trace zero. Assuming that there is a
nonegative integer} $m$ s\textit{uch that }$C_{m}$ \textit{or} $D_{m}$ 
\textit{is zero, then there is a positive integer} $k$ \textit{such that} $%
z^{k}=1$ \textit{and} $z$ \textit{is} $\left( k+1\right) $\textit{-potent
element.\smallskip }

\textbf{Proof. }Since $w$ has the trace zero, let $w^{2}=\beta $, with $\tau 
$ the multiplicative order of $\beta $. We have that $%
z^{m}=C_{m}+D_{m}w,C_{m},D_{m}\in \mathbb{Z}_{p}$. Supposing that $C_{m}$ is
zero, then we have $z^{m}=D_{m}w,$ with $\theta $ the multiplicative order
of $D_{m}$. Therefore $z^{2m}\in \mathbb{Z}_{p}^{\ast }$, and let $%
M=o(z^{2m})$, the order of the element $z^{2m}$. Therefore, $k=2mM$, if $%
z^{2m}\neq 1$ and $k=2m$, if $z^{2m}=1$. If $D_{m}$ is zero, then we have $%
z^{m}=C_{m}$ with $\upsilon $ the multiplicative order of $C_{m}$. It
results that $z^{\upsilon m}=1$.\smallskip

Now, we can say that we proved the following theorem.\smallskip

\textbf{Theorem 11.} \textit{With the above notations, an element} $z\in
A_{t}$ \textit{is a} $k$\textit{-potent element, if} $z$ \textit{is of one
of the forms:}

\textit{Case 1.} $n_{z}\neq 0$.

\textit{i)} $z=1+w$\textit{, with} $w\in A_{t}$\textit{,} $w$ \textit{is a
nilpotent element. In this case, }$z$ \textit{is} $\left( p+1\right) $%
\textit{-potent;}

\textit{ii)} $z=1+w$\textit{, with} $w\in A_{t}$ \textit{sucht that} $%
w^{2}=-1$\textit{. Since} $z^{4}=-4$ \textit{modulo} $p$ \textit{and} $%
\theta $ \textit{is the multiplicative order of }$-4$ \textit{in} $\mathbb{Z}%
_{p}^{\ast }$\textit{, we have that }$z$\textit{\ is }$\left( 4\theta
+1\right) $\textit{-potent}.

\textit{iii)} $z=a+w$\textbf{, }\textit{where }$a\in \mathbb{Z}_{p}$\textit{%
, }$w\in A_{t}$\textit{\ with }$t_{w}=0$\textit{,} $w^{2}=\beta \in \mathbb{Z%
}_{p}^{\ast }$\textit{, with} $\tau $ \textit{the multiplicative order of} $%
\beta $, \textit{and} $z^{r}=C_{r}+D_{r}w$.\textit{\ Assuming that there is
a nonegative integer} $m$ s\textit{uch that }$C_{m}$ \textit{or} $D_{m}$ 
\textit{is zero, then there is a positive integer} $s$ \textit{such that} $%
z^{s}=1$ \textit{and} $z$ \textit{is} $\left( s+1\right) $\textit{-potent
element.\smallskip\ If }$C_{m}=0$\textit{, then} $s=2mM$\textit{, where} $%
M=o(z^{2m})$, \textit{the order of the element} $z^{2m},$ \textit{if} $%
z^{2m}\neq 1$ \textit{and} $s=2m$, \textit{if} $z^{2m}=1$. \textit{\ If} $%
D_{m}=0$\textit{, then we have }$s=\upsilon m$\textit{, with} $\upsilon $ 
\textit{the multiplicative order of} $C_{m}$\textit{. Then }$k=s+1.$

\textit{Case 2.} $n_{z}=0$\textit{. The element} $z\in A_{t}$ \textit{is} $k$%
\textit{-potent if and only if} $t_{z}$ \textit{is} $k$\textit{-potent
element in} $\mathbb{Z}_{p}^{\ast }$\textit{, that means} $k-1$ \textit{is a
divisor of} $p-1$.\smallskip

\textbf{Example 14. }In the following, we will give some examples of values
of the potency index $k$.

i) Case $p=5$ and $t=2$, therefore we work on quaternions. We consider $%
z=2+i+j+k$ with the norm $n_{x}=2\neq 0$. We have $w=i+j+k$ and $%
z=2+w,w^{2}=2$. We have $z^{2}=1+4w,z^{3}=4w$, therefore $z^{6}=2,$ with $%
o\left( z^{6}\right) =o\left( 2\right) =4=M$. Since $m=3$, we have that $%
z^{24}=1$, then \thinspace $z^{25}=z$ and $z$ is $25$-potent element, $s=24,$
$k=25$.

ii) Case $p=7$, $t=2$ and $z=2+i+j+k$. The norm is zero and the trace is $4$%
. Since $4$ has multiplicative order equal with $3$, from Proposition 4, we
have $z^{4}=z$. Indeed, $z^{2}=1+4w,z^{3}=4+2w,z^{4}=2+w=z$ and $k=4$.

iii) Case $p=5$ and $t=2$. The element $z=1+3i+4j$ has $n_{z}=1,w=3i+4j$,
with $n_{w}=t_{w}=0$, therefore \thinspace $w$ is a nilpotent
element\thinspace . We have $z^{5}=1$, $z^{6}=z$ and $s=5,$ $k=6$.

iv) Case $p=3$ and $t=2$. The element $z=1+i+j+k$ has $n_{z}=1$ and $w=i+j+k$%
. We have $z^{2}=\left( 1+w\right) ^{2}=1+2w$, $z^{3}=\left( 1+w\right)
\left( 1+2w\right) =1+2w+w=1$, therefore $z^{4}=z$ and $s=3,k=4.$

v) Case\textbf{\ }$p=5,t=2$. We consider the element\textbf{\ }$%
z=2+3i+j+3k=2+3w,w=i+2j+k,n_{z}=3,n_{w}=1,t_{w}=0$\textbf{, }then\textbf{\ }$%
w^{2}=-1$\textbf{. }We have that $z^{2}=2w$\textbf{. }Therefore\textbf{\ }$%
z^{4}=1.$\textbf{\ }It results\textbf{\ }$k=4$\textbf{. }

vi) Case $p=5,t=2$. We consider the element $z=2+i+j+k=2+w$ with $%
n_{z}=2,n_{w}=3,t_{w}=0,w^{2}=2$ and $\tau =4$, the order of $\beta =2$. We
have $z^{2}=1+4w,z^{3}=4w,z^{6}=2,$ then and $M=4$. It results that $%
s=24,z^{24}=1$, then $z^{25}=z$ and $k=25$.

vii) Case $p=11,t=2$. We consider the element $z=2i+7j+3k$ with $%
n_{z}=7,z^{2}=4$, therefore $m=2,D_{2}=0,C_{2}=4,\upsilon =5$, the
multiplicative order of $C_{2}=4$. We have $z^{m\upsilon }=z^{10}=1$ and $%
k=11$.

viii) Case $p=13,t=3$, therefore we work on octonions. We consider the
element $z=3+2f_{1}+f_{2}+f_{3}+f_{4}+f_{5}+f_{6}+f_{7}=3+w$, $%
w=2f_{1}+f_{2}+f_{3}+f_{4}+f_{5}+f_{6}+f_{7}$, with $n_{z}=6$, $%
n_{w}=10,t_{w}=0$. We have $w^{2}=3$ and $\tau =3$, the order of $\beta =3$.
It results, $z^{2}=12+6w,z^{3}=2+4w,z^{6}=\left( 2+4w\right) ^{2}=3w$, $m=6$%
, and $z^{12}=1$. Therefore \thinspace $s=12$. We get $z^{13}=z$ and $k=13$.

ix) Case $p=17,t=4$, therefore we work on sedenions. The Sedenion algebra is
a noncommutative, nonassociative and nonalternative algebra of dimension $16$%
. We consider the element $z=1+w,w=\underset{i=1}{\overset{15}{\sum }}f_{i}$%
, with $w^{2}=2$ and $\tau =8$. It results $z^{2}=3+2w,z^{3}=7+5w,z^{4}=12w$%
. Then $m=4$, $z^{8}=16=2^{4}$ and $M=o\left( 16\right) =2$. It results
\thinspace $s=16$ and $k=17$.\smallskip

\textbf{Remark 15.} The $\left( m,k\right) $-type elements in $A_{t}$, with $%
m$, $n$ positive integers, are the elements $x\in A_{t}$ such that $%
x^{m}=x^{k}$, $m\geq k,$ smallests with this property. If $n_{x}\neq 0$,
then $x^{m-k}=1$ and $x$ is an $\left( m-k+1\right) $-potent element. If $%
n_{x}=0$ and $t_{x}\neq 0$, we have that $t_{x}^{m-k}=1$, then $x$ is an $%
\left( m-k+1\right) $-potent element. Therefore, an $\left( m,k\right) $%
-type element in $A_{t}$ is an $\left( m-k+1\right) $-potent element in $%
A_{t}$.

\begin{equation*}
\end{equation*}

\textbf{4. A nonunitary ring structure of quaternion Fibonacci elements over 
}$\mathbb{Z}_{p}$%
\begin{equation*}
\end{equation*}

The Fibonacci \ numbers was introduced \ by \textit{Leonardo of Pisa
(1170-1240) }in his book \textit{Liber abbaci}, book published in 1202 AD
(see [Kos; 01], p. 1-3). The $n$th term of these numbers is given by the
formula:%
\begin{equation*}
f_{n}=f_{n-1}+f_{n-2,}\ n\geq 2,\ 
\end{equation*}%
where $f_{0}=0,f_{1}=1.\medskip $

In [Ho; 63], were defined \ and studied Fibonacci quaternions over $\mathbb{H%
}$, defined as follows%
\begin{equation*}
F_{n}=f_{n}1+f_{n+1}i+f_{n+2}j+f_{n+3}k,
\end{equation*}%
called the \ $n$th Fibonacci quaternions.

In the same paper, the norm formula \ for the \ $n$th Fibonacci quaternions
was found:

\begin{equation*}
\boldsymbol{n}\left( F_{n}\right) =F_{n}\overline{F}_{n}=3f_{2n+3},
\end{equation*}%
where \ $\overline{F}_{n}=f_{n}\cdot 1-f_{n+1}i-f_{n+2}j-f_{n+3}k$ is the
conjugate of the $F_{n}$ in the algebra $\mathbb{H}.$

Fibonacci sequence is also studied when it is reduced modulo $m$. This
sequence is periodic and this period is called \textit{Pisano's period}, $%
\pi \left( m\right) $. In the following, we consider $m=p$, a prime number
and $\left( f_{n}\right) _{n\geq 0},$ the Fibonacci numbers over $\mathbb{Z}%
_{p}$. It is clear that, in general, the sum of two arbitrary Fibonacci
numbers is not a Fibonacci numbers, but if these numbers are consecutive
Fibonacci numbers, the sentence is true. In the following, we will find
conditions when the product of two Fibonacci numbers is also a Fibonacci
number.\smallskip\ In the following, we work on $A_{t},t=2$, over the field $%
\mathbb{Z}_{p}$. We denote this algebra with $\mathbb{H}_{p}$.

Let $F_{1}=a+bi+\left( a+b\right) j+\left( a+2b\right) k$ and $%
F_{2}=c+di+\left( c+d\right) j+\left( c+2d\right) k$, two Fibonacci
quaternions in \thinspace $\mathbb{H}_{p}$. \ We will find conditions such
that $F_{1}F_{2}$ and $F_{2}F_{1}$ are also Fibonacci quaternion elements,
that means elements of the same form: 
\begin{equation}
A+Bi+\left( A+B\right) j+\left( A+2B\right) k.  \tag{10.}
\end{equation}

We compute $F_{1}F_{2}$ and $F_{2}F_{1}$ and we obtain that 
\begin{equation}
F_{1}F_{2}=\left( -ac-3ad-3bc-6bd\right) +2adi+2a\left( c+d\right) j+\left(
2ac+ad+3bc\right) k  \tag{11.}
\end{equation}%
and 
\begin{equation}
F_{2}F_{1}=\left( -ac-3ad-3bc-6bd\right) +2bci+2c\left( a+b\right) j+\left(
2ac+3ad+bc\right) k.  \tag{12.}
\end{equation}%
By using relation $\left( 10\right) $, we get the following systems, with $%
c,d$ as unknowns. From relation $\left( 11\right) $, we obtain:

\begin{equation}
\left\{ 
\begin{array}{c}
\left( -3a-3b\right) c+\left( -3a-6b\right) d=0 \\ 
\left( -6b-3a\right) c+\left( -6b\right) d=0%
\end{array}%
\right.  \tag{13.}
\end{equation}

From relation $\left( 12\right) $, we obtain the system:%
\begin{equation}
\left\{ 
\begin{array}{c}
\left( -3a+3b\right) c+\left( -3a\right) d=0 \\ 
\left( -3a\right) c+\left( -6a-6b\right) d=0%
\end{array}%
\right.  \tag{14.}
\end{equation}

We remark that for $p=3$, the systems $\left( 13\right) $ and $\left(
14\right) $ have solutions, therefore, for $p=3$, there is a chance to
obtain an algebraic structure on the set $\mathcal{F}_{\pi \left( p\right) }$%
, the set of Fibonacci quaternions over $\mathbb{Z}_{p}.$

For $p=3$, the Pisano's period is $8$, then we have the following Fibonacci
numbers: $0,1,1,2,0,2,2,1$. We obtain the following Fibonacci quaternion
elements: $F_{0}=i+j+2k,F_{1}=1+i+2j,F_{2}=1+2i+2k,F_{3}=2+2j+2k,$\newline
$F_{4}=2i+2j+k,F_{5}=2+2i+j,F_{6}=2+i+k,F_{7}=1+j+k$, therefore $\mathcal{F}%
_{\pi \left( p\right) }=\{F_{i},i\in \{0,1,2,3,4,5,6,7\}\}$. All these
elements are zero norm elements. $F_{0}$ and $F_{4}$ are nilpotents, $F_{3}$%
, $F_{5}$ and $F_{6}$ are idempotent elements, $F_{1},F_{2},F_{7}$ are $3$%
-potent elements, By usyng $C++$ software, we computed the sum and the
product of these $8$ elements. Therefore, we have $F_{0}F_{i}=0$, for $i\in
\{0,1,...,7\},F_{4}F_{i}=0$, for $i\in \{0,1,...,7\},F_{5}F_{i}=F_{i}$, for $%
i\in \{0,1,...,7\},F_{6}F_{i}=F_{i}$, for $i\in \{0,1,...,7\}$ and%
\begin{eqnarray*}
F_{1}F_{0} &=&F_{4},F_{1}^{2}=F_{5},F_{1}F_{2}=F_{6},F_{1}F_{3}=F_{7}, \\
F_{1}F_{4} &=&F_{0},F_{1}F_{5}=F_{1},F_{1}F_{6}=F_{2},F_{1}F_{7}=F_{3},
\end{eqnarray*}%
\begin{eqnarray*}
F_{2}F_{0} &=&F_{4},F_{2}F_{1}=F_{5},F_{2}^{2}=F_{6},F_{2}F_{3}=F_{7}, \\
F_{2}F_{4} &=&F_{0},F_{2}F_{5}=F_{1},F_{2}F_{6}=F_{2},F_{2}F_{7}=F_{3,}
\end{eqnarray*}%
\begin{eqnarray*}
F_{3}F_{0} &=&F_{0},F_{3}F_{1}=F_{1},F_{3}F_{2}=F_{2},F_{3}^{2}=F_{3}, \\
F_{3}F_{4} &=&F_{4},F_{3}F_{5}=F_{5},F_{3}F_{6}=F_{6},F_{3}F_{7}=F_{7},
\end{eqnarray*}%
\begin{eqnarray*}
F_{7}F_{0} &=&F_{4},F_{7}F_{1}=F_{5},F_{7}F_{2}=F_{6},F_{7}F_{3}=F_{7}, \\
F_{7}F_{4} &=&F_{0},F_{7}F_{5}=F_{1},F_{7}F_{6}=F_{2},F_{7}^{2}=F_{3}.
\end{eqnarray*}%
Regarding the sum of two Fibonacci quaternions over $\mathbb{Z}_{3}$, we
obtain: 
\begin{equation*}
2F_{0}=F_{4},F_{0}+F_{1}=F_{2},F_{0}+F_{2}=F_{7},F_{0}+F_{3}=F_{6},F_{0}+F_{4}=0,
\end{equation*}%
\begin{eqnarray*}
F_{0}+F_{5}
&=&F_{3},F_{0}+F_{6}=F_{5},F_{0}+F_{7}=F_{1},2F_{1}=F_{5},F_{1}+F_{2}=F_{3},
\\
F_{1}+F_{3}
&=&F_{0},F_{1}+F_{4}=F_{7},F_{1}+F_{5}=0,F_{1}+F_{6}=F_{4},F_{1}+F_{7}=F_{6},
\end{eqnarray*}%
\begin{eqnarray*}
2F_{2}
&=&F_{6},F_{2}+F_{3}=F_{4},F_{2}+F_{4}=F_{1},F_{2}+F_{5}=F_{0},F_{2}+F_{6}=0,
\\
F_{2}+F_{7}
&=&F_{5},2F_{3}=F_{7},F_{3}+F_{4}=F_{5},F_{3}+F_{5}=F_{2},F_{3}+F_{6}=F_{1},
\end{eqnarray*}%
\begin{eqnarray*}
F_{3}+F_{7}
&=&0,2F_{4}=F_{0},F_{4}+F_{5}=F_{6},F_{4}+F_{6}=F_{0},F_{4}+F_{7}=F_{2}, \\
2F_{5}
&=&F_{1},F_{5}+F_{6}=F_{7},F_{5}+F_{7}=F_{4},2F_{6}=F_{2},F_{6}+F_{7}=F_{0},
\\
2F_{7} &=&F_{3}.
\end{eqnarray*}

From here, we have the following result.\smallskip .

\textbf{Proposition 16.} \ $\left( \mathcal{F}_{\pi \left( 3\right) }\cup
\{0\},+\right) $ \textit{is an abelian group of order} $9$\textit{,
isomorphic to} $\mathbb{Z}_{3}\times \mathbb{Z}_{3}$ \textit{and} $\left( 
\mathcal{F}_{\pi \left( 3\right) }\cup \{0\},+,\cdot \right) $ \textit{is a
nonunitary and noncommutative ring.}%
\begin{equation*}
\end{equation*}

\textbf{5. An application in Cryptography}

\begin{equation*}
\end{equation*}

We consider an algebra $A_{t}$ over $\mathbb{Z}_{p}$. This algebra is of
dimension $2^{t}$. We suppose that we have a text $m$ to be encrypted and
the alphabet has $p$ elements, \ $p$ a prime number. To each letter from
alphabet, will corespond a label from $0$ to $p-1$, that means we work on $%
\mathbb{Z}_{p}$. The encryption algorithm is the following.

1) We will split $m$ in blocks and we will choose the lenght of the blocks
of the form $2^{t}$. For a fixed $t$, we will find an invertible element $%
q,q\in $\thinspace $A_{t}$, that means $n_{q}\neq 0$. This element will be
the encryption key.

2) Supposing that $m=m_{1}m_{2}...m_{r}$ is the plain text, with $m_{i}$
blocks of lenght $2^{t}$, formed by the labels of the letters, to each $%
m_{i}=m_{i0}m_{i1}...m_{i2^{t}-1}$ we will associate an element $v_{i}\in
A_{t},v_{i}=\overset{2^{t}-1}{\underset{j=0}{\sum }}m_{ij}f_{j}$.

3) We compute $qv_{i}=w_{i}$, for all $i\in \{1,2,...,r\}$. We obtain $%
w=w_{1}w_{2}...w_{r}$, the encrypted text.

To decrypt the text, we use the decryption key, then we compute $d=q^{-1}$
and $v_{i}=dw_{i}$, for all $i\in \{1,2,...,r\}$.\smallskip

\textbf{Example 17.} We consider the word MATHEMATICS and the key SINE. We
work on an alphabet with $29$ letters, including blank space, denoted with
"*", "." and ",". The labels of the letters are done in the below table

\begin{equation*}
\begin{tabular}{|l|l|l|l|l|l|l|l|l|l|}
\hline
A & B & C & D & E & F & G & H & I & J \\ \hline
0 & 1 & 2 & 3 & 4 & 5 & 6 & 7 & 8 & 9 \\ \hline
K & L & M & N & O & P & Q & R & S & T \\ \hline
10 & 11 & 12 & 13 & 14 & 15 & 16 & 17 & 18 & 19 \\ \hline
U & V & W & X & Y & Z & * & . & , &  \\ \hline
20 & 21 & 22 & 23 & 24 & 25 & 26 & 27 & 28 &  \\ \hline
\end{tabular}%
\end{equation*}

We consider $t=2$, therefore we work on quaternions. We will add an "A" at
the end of word "MATHEMATICS", to have multiple of $4$ lenght text,
therefore, we will encode the text "MATHEMATICSA". We have the following
blocks MATH, EMAT, ICSA, with the corresponding quaternions $v_{1}=12+19j+7k$%
, for MATH, $v_{2}=4+12i+19k$, for EMAT and $v_{3}=8+2i+18j$ for ICSA. The
key is $q=18+8i+13j+4k$, it is an invertible element, with the nonzero norm, 
$n_{q}=22$. We have $w_{1}=qv_{1}=28+24i+7j+7k$, corresponding to the
message ",YHH", $w_{2}=qv_{2}=16+2i+6j+28k$, corresponding to the message
"QCG," and $w_{3}=qv_{3}=10+28i+j+5k$, corresponding to the message "K,BF".
Therefore, the encrypted message is ",YHHQCG,K,BF". The decryption key is $%
d=q^{-1}=14+26i+6j+13k$. For decryption, we will compute $%
dw_{1}=12+19j+7k=v_{1}$, $dw_{2}=4+12i+19k=v_{2}$, $dw_{3}=8+2i+18j=v_{3}$,
and we find the initial text "MATHEMATICSA".

\begin{equation*}
\end{equation*}

\textbf{Conclusion.} In this paper we studied properties of some special
elements in algebras obtained by the Cayley-Dickson process and we find an
algebraic structure(nonunitary and noncommutative ring) over Fibonacci
quaternions over $\mathbb{Z}_{3}$. Moreover, an encryption method by using
these elements is also provided. As a further research, we intend to study
other special elements in the idea of {}{}finding another good properties. 
\begin{equation*}
\end{equation*}

\textbf{References}%
\begin{equation*}
\end{equation*}

\textbf{[}Ba; 01\textbf{]} Baez, J.C., \textit{The Octonions}, B. Am. Math.
Soc., \textbf{39(2)}(2001), 145-205, http://www.ams.org/journals/bull/

2002-39-02/S0273-0979-01-00934-X/S0273-0979-01-00934-X.pdf\smallskip .

[Ho; 63] Horadam, A. F., \textit{Complex Fibonacci Numbers and Fibonacci
Quaternions}, Amer. Math. Monthly \textbf{70}(1963), 289-291.\smallskip

[Kos; 01] ~Koshy, T., \textit{Fibonacci and Lucas Numbers with Applications}%
, A Wiley-Interscience publication, U.S.A, 2001.\smallskip

[MS; 11] Miguel, C. J., Serodio R., \textit{On the Structure of Quaternion
Rings over} $\mathbb{Z}p$, International Journal of Algebra, 5(27)(2011),
1313-1325.

\textbf{[}McC; 80\textbf{]} McCrimmon, K., \ \textit{Pre-book on Alternative
Algebras}, 1980,\newline
http://mysite.science.uottawa.ca/neher/Papers/alternative/\newline
http://mysite.science.uottawa.ca/neher/Papers/alternative/\newline
2.2.Composition\%20algebras.pdf.\smallskip

\textbf{[}Mo; 15\textbf{]} Mosi\'{c} Dijana, \textit{Characterizations of
k-potent elements in rings}, Annali di Matematica, 194(2015), 1157--1168,
DOI 10.1007/s10231-014-0415-5

\textbf{[}RPC; 22\textbf{]} Ratanaburee, P., Petapirak, M., Chuysurichay,
S., \textit{On (m, k) -type elements in the ring of integers modulo n},
Songklanakarin J. Sci. Technol., 44 (5)(2022), 1179--1184.

\textbf{[}Sc; 66\textbf{]} Schafer, R. D., \textit{An Introduction to
Nonassociative Algebras,} Academic Press, New-York, 1966.\smallskip

\textbf{[}Sc; 54\textbf{]} Schafer, R. D., \textit{On the algebras formed by
the Cayley-Dickson process,} Amer. J. Math., \textbf{76}(1954),
435-446.\smallskip

\textbf{[}Wu; 10\textbf{]} Y. Wu, \textit{k-Potent Matrices - Construction
and Applications in Digital Image Encryption}, Recent Advances in Applied
Mathematics, Proceedings of the 2010 American Conference on Applied
Mathematics, USA, 2010, 455--460.

\begin{equation*}
\end{equation*}

Cristina FLAUT

Faculty of Mathematics and Computer Science, Ovidius University,

Bd. Mamaia 124, 900527, Constan\c{t}a, Rom\^{a}nia,

http://www.univ-ovidius.ro/math/

e-mail: cflaut@univ-ovidius.ro; cristina\_flaut@yahoo.com

\begin{equation*}
\end{equation*}

Andreea BAIAS

PhD student at Doctoral School of Mathematics,

Ovidius University of Constan\c{t}a, Rom\^{a}nia,

e-mail: andreeatugui@yahoo.com

\end{document}